\documentclass[12pt,a4paper]{article}
\usepackage[english]{babel}
\usepackage[utf8x]{inputenc}
\PrerenderUnicode{Ž}
\PrerenderUnicode{ä}
\usepackage{amsmath}
\usepackage[margin=1in]{geometry}
\usepackage{graphicx}
\usepackage[section]{placeins}
\usepackage{float}
\usepackage{ amssymb }
\usepackage{hyperref}
\usepackage{amsthm}
\usepackage{amsfonts}
\usepackage{epsfig}
\usepackage{ textcomp }
\usepackage{mathtools}
\usepackage{authblk}
\usepackage{etoolbox}

\newcommand{\bb}{\begin{equation}}
\newcommand{\ee}{\end{equation}}

 \newtheorem{thm}{Theorem}
 
 \newtheorem{lem}[thm]{Lemma}
 \newtheorem{conj}[thm]{Conjecture}

\newcommand{\QED} {\hfill$\square$}

\DeclareMathOperator{\conv}{conv}

\title{Inequalities on Projected Volumes}
\date{}
\author{Imre Leader, Žarko Ranđelović, Eero Räty
}
\affil{\small{ Centre for Mathematical Sciences
\\
Wilberforce Road\\
Cambridge CB3 0WB, U.K.}}
\affil{\texttt{i.leader@dpmms.cam.ac.uk, zr233@cam.ac.uk, epjr2@cam.ac.uk}}

\begin{document}

\maketitle

\begin{abstract}
In this paper we study the following geometric problem: given $2^n-1$ real numbers $x_A$ indexed by the non-empty subsets $A\subset \{1,..,n\}$, is it possible to construct a body $T\subset \mathbb{R}^n$ such that $x_A=|T_A|$ where $|T_A|$ is the $|A|$-dimensional volume of the projection of $T$ onto the subspace spanned by the axes in $A$? As it is more convenient to take logarithms we denote by $\psi_n$ the set of all vectors $x$ for which there is a body $T$ such that $x_A=\log |T_A|$ for all $A$. Bollobás and Thomason showed that $\psi_n$ is contained in the polyhedral cone defined by the class of ‘uniform cover inequalities'. Tan and Zeng conjectured that the convex hull $\conv (\psi_n)$ is equal to the cone given by the uniform cover inequalities.\\
\\
We prove that this conjecture is `nearly' right: the closed convex hull $\overline{\conv }(\psi_n)$ is equal to the cone given by the uniform cover inequalities. However, perhaps surprisingly,  we also show that $\conv (\psi_n)$ is not closed for $n\ge 4$, thus disproving the conjecture.
\end{abstract}
\section{Introduction} 
Let $T$ be a body in $\mathbb{R}^n$, meaning a compact subset of $\mathbb{R}^n$. 
Let $\{e_1,..,e_n\}$ be the standard basis of $\mathbb{R}^n$ and let Span($A$) be the subspace spanned by $\{e_i$ $|i\in A\}$. Given a non-empty set $A\subset [n]=\{1,2,..,n\}$ with $|A|=d$, we denote by $T_A$ the projection of $T$ onto Span($A$) and let $|T_A|$ be its $d$-dimensional volume. We define $x(T)$ to be the log projection vector for a body $T$, meaning the $(2^n-1)$-dimensional vector with entries indexed by non-empty $A\subset [n]$ and $x(T)_A=\log |T_A|$. Note that if a body $T$ has positive volume then also $|T_A|>0$ for any $A$. Whenever we mention a $(2^n-1)$-dimensional vector we shall assume that the coordinates are indexed by the non-empty subsets of $[n]$. Say that a $(2^n-1)$-dimensional vector $x$ is \textit{constructible} if there is a body $T$ such that $x=x(T)$. The \textit{constructible region} $\psi_n$ is the set of all constructible vectors:
$$\psi_n=\{x\in\mathbb{R}^{2^n-1}\ |x\text{ is constructible}\}$$
What is the structure of $\psi_n$? As we will see, $\psi_n$ is not convex for $n\ge 4$. (This was stated by Tan and Zeng [6], but their proof is incorrect.) Bollobás and Thomason [2] proved a class of linear inequalities that hold for all constructible vectors. To state their theorem we first need a definition. We say that a collection $Y_1,..,Y_l$ of subsets of $[n]$ (not necessarily distinct) is a \textit{$k$-uniform cover} of $[n]$ if every element in $[n]$ appears in precisely $k$ of the $Y_i$. A cover is \textit{uniform} if it is $k$-uniform for some $k$.
 
 \begin{thm}
 (Bollobás-Thomason uniform cover theorem.) Let $T$ be a body and $Y_1,..,Y_l$ be a $k$-uniform cover of $[n]$. Then $$\prod\limits_{i=1}^l |T_{Y_i}|\ge |T|^k$$
 \end{thm}
 The discrete version of this theorem with set systems was first proved by Shearer -- see Chung, Frankl, Graham and Shearer [3].
 \begin{thm}
 (Product theorem.) Let $A_1,..,A_m$ be subsets of $[n]$ such that every element of $[n]$ is contained in at least $k$ of $A_1,..,A_m$. Let $\mathcal{F}$ be a collection of subsets of $[n]$ and let $\mathcal{F}_i=\{F\cap A_i:F\in \mathcal{F}\}$ for $1\le i\le m$. Then 
 $$ \prod\limits_{i=1}^m |\mathcal{F}_i|\ge |\mathcal{F}|^k$$
 \end{thm}
 If we take logarithms in Theorem 1 we get a set of linear inequalities that hold for all vectors in $\psi_n$. Bollobás and Thomason also showed that there are in fact only finitely many irreducible uniform covers (meaning the ones that cannot be decomposed into two smaller uniform covers), and of course any uniform cover has a decomposition into finitely many irreducible covers. This means that above we have in fact only a finite set of inequalities. \\
 \\
 Can these inequalities determine $\overline{\conv}(\psi_n)$? Not quite, as any positive linear combination (when taking logs) of the above inequalities will have only $\log |T|$ on the right-hand side but there are more inequalities that have to hold, such as $|T_1||T_2|\ge |T_{12}|$. This holds since if $T$ is a body then so is $T_{12}$, where $12$ in the index abbreviates $\{1,2\}$. We thus get more inequalities just by considering the uniform cover inequalities for all subsets of $[n]$. Defining in a similar way as above a $k$-uniform cover of an arbitrary $Y$ which is a subset of $[n]$, we get that if $Y_1,..,Y_l$ is a $k$-uniform cover of $Y$ then
 $$\prod_{i=1}^l |T_{Y_i}|\ge |T_Y|^k$$
 These are all still defined by a finite set of inequalities, because for each of the finitely many $Y$ there are only finitely many irreducible covers of $Y$. One can easily check that these are the only possible linear homogeneous inequalities with only one term on the right-hand side. Now, there are many other guesses of plausible inequalities, such as perhaps  $|T_{12}||T_{23}||T_{34}|\ge |T_{123}||T_{234}|$. However, after a bit of checking they all seem to fail.  In view of the above we define the set $BT_n$ as follows:
 \begin{align*}
 BT_n=\{x| \textrm{ for all }Y\subset [n],k\in \mathbb{N},Y_1,..,Y_l\textrm{ a }k\textrm{-uniform cover of Y,\textrm{ we have } }\sum_ix_{Y_i} \ge kx_Y \}
 \end{align*}
 Tan and Zeng [6] conjectured the following:
 \begin{conj}
 \emph(\emph[\emph6\emph]\emph).
 The set $\conv (\psi_n)$ is equal to $BT_n$
 \end{conj}
 In this paper we show that their conjecture is nearly correct. This will be our main result and it is the following theorem:
 \begin{thm}
 The set $\overline{\conv}(\psi_n)$ is equal to $BT_n$.
 \end{thm}
 This theorem says that in fact the uniform cover inequalities are the best we can do. It also provides us with only finitely many inequalities to check if we want to determine whether or not a given vector lies in $\overline{\conv}(\psi_n)$ (i.e $\overline{\conv}(\psi_n)$ is a polyhedral cone). \\
 \\
 The plan of this paper is as follows. In Section 2 we prove Theorem 4. In Section 3 we will prove various properties of $\psi_n$. We show that $\psi_n$ is not convex for $n\ge 4$. This was stated by Tan and Zeng [6], but their proof is incorrect. We also show that if $n\ge 4$ then $\conv(\psi_n)$ is, perhaps surprisingly, not closed. Further, we show a scaling property of $\psi_n$ which gives an unexpected second proof of Theorem 4.
 
  \section{Proof of Theorem 4}

 We note that $\overline{\conv}(\psi_n)\subset BT_n$ by the uniform cover inequalities. Thus we need to show that
$BT_n \subset \overline{\conv}(\psi_n)$. Define a \textit{box} in $\mathbb{R}^n$ to be a body of the form $\prod\limits_{i=1}^n [a_i,b_i]$, where $a_i<b_i$ for all $i$. We first need a simple lemma.
 \begin{lem}
 Suppose that $\alpha ,\beta_1,..,\beta_k$ are positive reals and consider the inequality $$\sum_{i=1}^k\beta_ix_{Y_i}\ge \alpha x_{[n]}$$ where $Y_i$ are subsets of $[n]$. If for each $1\le s\le n$ we have $\sum_{i;s\in Y_i}\beta_i=\alpha$ then the above inequality is a positive linear combination of uniform cover inequalities.
 \end{lem}
 \textit{Proof.} If all the $\beta_i$ are rational then so is $\alpha$ and multiplying by an integer gives a uniform cover inequality. Consider the following equations in $\mathbb{R}^k$:  $$\textrm{ each}\  2\le s\le n \sum_{i:1\in Y_i}v_i=\sum_{i:s\in Y_i}v_i$$ The $\beta_i$ are a solution to this set of equations. These equations all have rational coefficients i.e they can be represented in the form $M\bold{v}=0$ for some $(n-1)$ x $k$ rational matrix. Since the entries of $M$ are all in $\mathbb{Q}$ the space of rational solutions has a basis in $\mathbb{Q}^k$. The dimension of the space of solutions is equal to the nullity of $M$, which is the same over $\mathbb{Q}$ and $\mathbb{R}$, and rational vectors are linearly independent over $\mathbb{Q}$ if and only if they are linearly independent over $\mathbb{R}$. This means that the space of solutions over $\mathbb{R}$ has a basis of rational vectors. So we can express $(\beta_i)_{i=1}^k$ in terms of those rational vectors.\\
 \\
 Let $\bold{f_1,..,f_l}$ be the rational basis of the space of solutions (now over $\mathbb{R}$) and let $\bold{b}$ be the vector with coordinates $\beta_i$ and $\bold{b}=\sum_{i=1}^l \lambda_i\bold{f_i}$. Now consider rationals $q_{ij}$ such that $q_{ij}\rightarrow \lambda_i$ as $j\rightarrow \infty$ for each $i$. Let $\beta_{ij}$ be the $i$-th coordinate of the vector $\sum_{p=1}^l q_{pj}\bold{f_p}$. Now we have that $\beta_{ij}\rightarrow \beta_i$ as $j\rightarrow \infty$ for each $i$, so we may assume that all $\beta_{ij}$ are positive. But we also have that all $\beta_{ij}$ are rational. So if $1\le s\le n$ we can define $\alpha_j=\sum_{i;s\in Y_i}\beta_{ij}$ which is well defined since the right-hand side is the same for all $s$ by choice of $\beta_{ij}$. Now the inequalities $$\sum_{i=1}^{k}\beta_{ij}x_{Y_i}\ge \alpha_jx_{[n]}$$ are all uniform cover inequalities (when scaled) since all the coefficients are rational but they also imply our original inequality. Thus Farkas' Lemma (see e.g.$\ $[5]) gives the desired result. 
\QED
\\

 \textit{Proof of Theorem 4.} Suppose we are given a linear homogeneous inequality which holds for all vectors in $\psi_n$ and hence also in $\overline{\conv}(\psi_n)$. We may write it in the form 
 \begin{align}
 \sum_{i=1}^l \alpha_ix_{A_i}\ge \sum_{i=1}^k \beta_ix_{B_i}
 \end{align}
 where all $\alpha_i$ and $\beta_i$ are positive. Now consider variables $b_{ij}$ for each $j\in B_i$ and $a_{ij}$ for each $j\in A_i$. Now consider the set of inequalities on the variables $a_{ij},b_{ij}$ given by each $i',i$ such that $A_{i'}\subset B_i$ where the inequality is 
 \begin{align} 
 \sum_{j\in A_{i'}}a_{i'j}\ge \sum_{j\in A_{i'}}b_{ij}
 \end{align}
 We will now split into two cases depending on whether or not the inequalities in (2) imply the following inequality:
 \begin{align} 
 \sum_{i,j\in A_{i}}\alpha_ia_{ij}\ge \sum_{i,j\in B_{i}}\beta_ib_{ij}
 \end{align}
 Supposing this was true. We get that by Farkas' lemma the inequality in $(3)$ is a positive linear combination of inequalities in $(2)$. But now we can see that $(3)$ is in fact a positive linear combination of uniform cover inequalities. This is because if we fix $i$ then the right-hand side of $(3)$ has the term $S_i=\beta_i\sum_{j\in B_i}b_{ij}$. The only inequalities from $(2)$ that can contribute to that term are those with $A_{i'},B_i$ for some $A_{i'}\subset B_i$. Every inequality in $(2)$ contributes to exactly one of the $S_i$. Now we fix $i$ and suppose that for each $A_j\subset B_i$ we have that $\lambda_j\ge 0$ is the coefficient of the inequality in $(2)$ with $A_j,B_i$ in the positive linear combination giving $(3)$. Then consider the following inequality on the vectors $x\in \mathbb{R}^{2^n-1}$ 
 \begin{align} 
 \sum_{j,A_j\subset B_i}\lambda_jx_{A_j}\ge \beta_ix_{B_i}
 \end{align} 
 For all $s\in B_i$ we have $\sum\limits_{j,s\in A_j}\lambda_j=\beta_i$ by looking at the coefficient of $b_{is}$. By Lemma $5$ this is a positive linear combination of uniform cover inequalities on the set $B_i$. This gives us the desired result since we can split the inequality in $(3)$ into $k$ inequalities, each corresponding to a set $B_i$ and this will give a splitting of $(1)$ into $k$ inequalities of the form in $(4)$ which are all positive linear combinations of uniform cover inequalities. \\
 \\
 Now suppose that the inequalities in $(2)$ do not imply the inequality in $(3)$. Consider some $b_{ij},a_{ij}$ which satisfy all the inequalities in (2) but such that $\sum_{i,j\in B_{i}}\beta_ib_{ij}>\sum_{i,j\in A_{i}}\alpha_ia_{ij}$. Let $M=\log{k}\sum_i \alpha_i+1$. By scaling we can assume that $\sum_{i,j\in B_{i}}\beta_ib_{ij}>M+\sum_{i,j\in A_{i}}\alpha_ia_{ij}$. Now consider the set $T$ which is the union of the boxes $X_1,..,X_k$ where $X_i$ is a $|B_i|$-dimensional box in Span$(B_i)$ with sidelength $e^{b_{ij}}$ in direction $j$ for each $j\in B_i$. Also let $x=x(T)$ where $x(T)$ is defined. Then $x_{B_i}\ge \sum_{j\in B_i}b_{ij}$ so the right-hand side of $(1)$ is at least $ \sum_{i=1}^k \beta_i\sum_{j\in B_i}b_{ij}$ For the left-hand side terms we have that $$x_{A_i}\le \log{k}+\max_{i';A_i\subset B_{i'}} \sum_{j\in A_i}b_{i'j}\le \log{k}+\sum_{j\in A_i}a_{ij}$$ The last inequality is due to $(2)$ and the first is from the fact that $X_j$ contributes to the volume of $T_{A_i}$ only if $A_i\subset B_j$. But now it is clear that $$\sum_{i=1}^l\alpha_ix_{A_i}< M+\sum_{i,j\in A_{i}}\alpha_ia_{ij}\le \sum_{i,j\in B_{i}}\beta_ib_{ij}\le \sum_{i=1}^k \beta_ix_{B_i} $$ which violates $(1)$. Since $T$ may not have positive volume we can add to it a tiny box of full dimension which still makes $(1)$ false and hence a contradiction.
 \\
 \\
 We have thus just shown that the set of linear homogeneous inequalities that hold for all points in $\overline{\conv}(\psi_n)$ is precisely that set of inequalities that hold for all points in $BT_n$. If there is an $x\in BT_n\setminus \overline{\conv}(\psi_n)$ then by the hyperplane separation theorem (see e.g.$\ $[1]) there is a hyperplane strongly separating $x$ and $\overline{\conv}(\psi_n)$ i.e there is a vector $w\in \mathbb{R}^{2^n-1}$ and $c\in \mathbb{R}$ such that for every $y\in \overline{\conv}(\psi_n)$ we have $w\cdot y< c$ and $w\cdot x>c$. To see this apply the regular hyperplane separation theorem  to $\overline{\conv}(\psi_n)$ and a small ball around $x$ disjoint from $\overline{\conv}(\psi_n)$. The box of full dimension with side lengths equal to $1$ gives that $\bold{0}\in \psi_n$ and hence $c>0$. Notice that for all bodies we have an inequality of the form $$c+\sum_{i=1}^l \alpha_ix_{A_i}\ge \sum_{i=1}^k \beta_ix_{B_i}$$ then by the same proof as above taking $M=\log{k}\sum_i \alpha_i+c+1$ we obtain that $$\sum_{i=1}^l \alpha_ix_{A_i}\ge \sum_{i=1}^k \beta_ix_{B_i}$$ is a positive linear combination of uniform cover inequalities. If we consider the inequality $w\cdot y\le c$ then we get from above that $w\cdot x\le 0$ since $x\in BT_n$ which is a contradiction.  \QED

 \section{Further properties of $\psi_n$}

 What else can we say about $\conv (\psi_n)$ or possibly $\psi_n$ itself? We will prove that for $n\ge 4$ $\conv (\psi_n)$ is not closed. The structure of $\psi_n$ itself is more complicated. We prove that $\psi_n$ is not convex and is not even a cone. Thus we have $$\psi_n\subsetneq\conv (\psi_n)\subsetneq \overline{\conv}(\psi_n)$$ We first prove that $\conv (\psi_n)$ is not closed (for $n\ge 4$). To do this we need a few technical lemmas to prove a result very similar to Theorem 4 in [2]. (The only difference with Theorem
4 in [2] is that here we are considering all compact sets, not just those that
are the closure of their interior.)
For a body $T$ define $T(x)$ to be the $(n-1)$-dimensional body which is the set of points in $T$ with last coordinate equal to $x$. We view $T(x)$ as a body in $\mathbb{R}^{n-1}$.
 \begin{lem}
 Let $T$ be a body and $a_k$ be a sequence of reals converging to some $a$. Suppose that $C>0$ and $|T(a_k)|\ge C$ for all $k$. Then $|T(a)|\ge C$.
 \end{lem}
 \textit{Proof.}
 Suppose this is not true. By regularity of measure there exists an open $U$ of Borel measure $m(U)<C$ containing $T(a)$. But then there are points $p_k\in T(a_k)$ not contained in $U$. By compactness $p_k$ has a subsequence converging to some $p$, but then by compactness $(p,a)\in T$ and hence $p\in U$, which is a contradiction. \QED
 \begin{lem}
 Let $A$ be a compact subset of $\mathbb{R}$. Then there exists a closed set $E\subset A$ such that $m(E)=m(A)$ and for every closed $B\subset A$ such that $m(B)=m(A)$ we have $E\subset B$.
 \end{lem}
 \textit{Proof.} Let $$E=\{x\in A|\textrm{ if }x\in U,\ U-\textrm{open then}\ m(U\cap A)\neq 0\}$$
 This set is trivially closed in $A$ and is hence closed in $\mathbb{R}$. If we list the rationals $q_1,q_2,...$ then for each $x\in A\backslash E$ there are some $q_i,q_j$ such that $x\in (q_i,q_j)$ and $m((q_i,q_j)\cap A)=0$. Summing over all countably many such $i,j$ we have that $m(A\backslash E)=0$. Suppose that $B$ is closed and is a subset of $A$. and let $x\in E\backslash B$. Then there is some open $U$, containing $x$, disjoint from $B$. But then $m(U\cap A)\neq 0$ and hence $m(B)\le m(U^c\cap A)<m(A)$. Thus if $m(B)=m(A)$ then $E\subset B$. \QED
\\
\\
 We remark that the same result holds for $\mathbb{R}^n$ and the proof is almost identical. \\
 \\
 Define an \textit{open product set} in $\mathbb{R}^n$ to mean a set $B$ of the form $B=\prod\limits_{i=1}^nB_i$ where the $B_i$ are bounded open subsets of $\mathbb{R}$. If the $B_i$ are compact call it a \textit{compact product set}.

 \begin{lem}
 Suppose that $A_1,..,A_s$ is a $k$-uniform cover of $[n]$ and let $E_i$ be equivalence classes under the relation of $i\sim j$ if for all $l$ either both or neither of $i$ and $j$ are in $A_l$. Now let $T$ be a body for which $\prod\limits_l |T_{A_l}|=|T|^k$. Then there are bodies $K_{E_i}\subset T_{E_i}$ such that $\prod\limits_i K_{E_i}\subset T$ and $|\prod\limits_i K_{E_i}|=|T|$.
\end{lem}
\textit{Proof.} Let us first assume for simplicity that all the equivalence classes are singletons $\{i\}$. We prove this statement by induction on $n$. The case $n=1$ is trivial. For $n\ge 2$ we shall call $T(x)$ to be the $(n-1)$-dimensional body consisting of all points in $T$ with $n$-th coordinate equal to $x$. Let $I=\{i| n\in A_i\}$ and $J=[s]\backslash I$. Then we have the following
$$|T(x)_{A_i}|\le |T_{A_i}|\ \ \textrm{for}\ i\in J$$
$$|T_{A_i}|=\int|T(x)_{A_i\backslash \{n\} }|dx\ \ \textrm{for}\ i\in I $$ The collection $A_i$ for $i\in J$ together with $A_i\backslash \{n\}$ for $i\in I$ form a $k$-uniform cover of $[n-1]$ whose equivalence classes are all singletons hence by using Holder's inequality and the uniform cover inequality we obtain

\begin{gather*} |T|=\int|T(x)|dx\le \int \bigg[ \prod\limits_{i\in I}|T(x)_{A_i\backslash \{n\}}| \prod\limits_{i\in J} |T(x)_{A_i}|\bigg]^{1/k}dx \\ \le \prod\limits_{i\in J}|T_{A_i}|^{1/k}\int\prod\limits_{i\in I} |T(x)_{A_i\backslash \{n\}}|^{1/k}dx\ \ \ \ \ \ \ \ \ \ \ \ \ \ \ \ \ \ \ \\ \le \prod\limits_{i\in J}|T_{A_i}|^{1/k}\bigg[\prod\limits_{i\in I}\int |T(x)_{A_i\backslash \{n\}}|dx\bigg]^{1/k}\ \ \ \ \ \ \ \ \ \ \ \ \ \ \ \ \ \\ \le \prod\limits_{i=1}^s|T_{A_i}|^{1/k}\ \ \ \ \ \ \ \ \ \ \ \ \ \ \ \ \ \ \ \ \ \ \ \ \ \ \ \ \ \ \ \ \ \ \ \ \ \ \ \ \ \ \ \ \ \ \ \ \ \ \ \end{gather*} The case $|T|=0$ is trivial so assume $|T|>0$. By the assumption we have equality everywhere above. Thus we know that \begin {align} |T(x)|=\prod\limits_{i\in I}|T(x)_{A_i\backslash \{n\}}|^{1/k} \prod\limits_{i\in J} |T(x)_{A_i}|^{1/k} \end{align} holds for almost all $x$. By induction for almost all $x$ we have that $T(x)$ contain a compact product set $M(x)=\prod_{i=1}^{n-1}M(x)_i$ with $|T(x)|=|M(x)|$. Using this and the second inequality we also have that $F'=\{x\in F|\textrm{ (5) holds and }|T(x)_{A_i}|=|T_{A_i}|\textrm{ for all }i\in J\}$ has the same measure as $F=\{x|\ |T(x)|>0\}$ Note that both sets are measurable because Lemma 6 gives that $\{x|\ |T(x)|\ge C\}$ is closed for any $C>0$ and with an identical proof so are $\{x|\ |T(x)_{A_i}|\ge C\}$. Now we may replace $T$ with $T\cap \mathbb{R}^{n-1}\times \overline{F'}$ which has the same volume as $T$ and hence also the same volumes of projections to coordinate sets $A_i$.\\
\\
Thus by Lemma 6 $|T(x)_{A_i}|=|T_{A_i}|$ for all $x\in E=\overline{F'}$ and all $i\in J$ and because of $(5)$ we also have $|M(x)_{A_i}|=|T(x)_{A_i}|$. Hence $$|M(x)_{A_i}|=|T_{A_i}|$$
But $M(x)_{A_i}$ are compact product sets and all contained in $T_{A_i}$. This means that for any $x,y\in E$ the `sides' of the compact product sets are essentially the same set i.e. $|M(x)_j|=|M(y)_j|=|M(x)_j\cap M(y)_j|$ for all $j\in A_i$ with $i\in J$. Since the equivalence classes are all singletons the $A_i$ with $i\in J$ cover $[n-1]$ and hence the above is true for all $j$. We are almost done. Now we pick an arbitrary $x_0\in E$ and let $K_i$ be the set given by applying Lemma 7 to $M(x_0)_i$ for $i\le n-1$. Also set $K_n=E$. By the above properties and since all $M(x)_i$ are compact it is easy to check that $\prod\limits_{i=1}^n K_i$ is a body contained in $T$ with the same volume as $T$. \\
\\
By the remark after Lemma 7 and a similar statement to Lemma 6 reducing multiple dimensions we deduce by the same proof as above (inducting on the number of classes) the case when not all classes are singletons. \QED
\\
\\
We remark that Lemma 8 is not true if we remove the assumption that $T$ is compact. Indeed, let $T=[0,1]^2\backslash \{(x,x)|x\in [0,1]\}$. If $A\times B\subset T$ then $A\cap B=\emptyset $ so $|A||B|\le \frac{1}{4}$.\newpage
\begin{thm}
If $n\ge 4$ then $\conv (\psi_n)$ is not closed.
\end{thm}
Consider the vector $v\in \mathbb{R}^{2^n-1}$ with $v_{24}=v_{13}=2$,\ $v_{123}=v_{234}=v_1=v_2=v_3=v_4=1$ and all other coordinates $0$. It is easy to show that all irreducible uniform covers hold on $v$. We also have equality in the uniform cover inequalities $2v_{123}=v_{12}+v_{13}+v_{23}$ and $2v_{234}=v_{23}+v_{24}+v_{34}$. If $v$ was in a convex combination of constructible vectors then all of those must satisfy equality in the above 2 uniform cover inequalities. Suppose that $w$ is one of those. Suppose that a body $T$ has log projection vector $w$. By Lemma $8$, $T_{123},T_{234}$ contain compact product sets of the same volume as those projections. The $23$-projections of those sets must have volume $|T_{23}|$ and hence must have equal measure in both the $x_2,x_3$ projections. Hence we must have, say, the equation $$w_{123}-w_{12}=w_{234}-w_{24}$$
This equation hence must hold for $v$ as well but this is not true. Thus $\conv (\psi_n)$ is not equal to its closure by Theorem $4$ and is hence not closed.\QED\\
\\
We now move on to results about $\psi_n$ itself. It turns out that the constructible region is not convex. For this we will need the following result of Ellis, Friedgut, Kindler and Yehudayoff (Corollary 2 in [4]). It is a stability version of Lemma 8.
\begin{lem}
\emph[\emph4\emph]
For every integer $n\ge 2$ there exists $b=b(n)>0$ such that the following holds. Let $m\in \mathbb{N}$ and let $\mathcal{F}\subset \mathcal{P}([n])$ be an $m$-uniform cover with $\sigma(\mathcal{F})>0$ where $\sigma(\mathcal{F})$ is the largest integer $\sigma$ such that for every $i,j\in [n]$ there are at least $\sigma$ sets containing $i$ but not $j$. Let $T$ be a bounded open set such that 
$$|T|\ge (1-\epsilon)\bigg(\prod\limits_{A\in \mathcal{F}}|T_A|\bigg)^{1/m}$$ Then there exists an open product set $B\subset \mathbb{R}^n$ such that $|T\triangle B|\le \frac{bm}{\sigma (\mathcal{F})}\epsilon |T|$
\end{lem}
The idea is that we will prove that there is a small neighbourhood around the above described $v$ which contains no constructible vectors. Since $\conv (\psi_n)$ is dense in $BT_n$ we will deduce that $\psi_n$ is not convex.  We remark that one can actually use Lemma 10 to give a proof of Lemma 8.

\begin{thm}
If $n\ge 4$ then $\psi_n$ is not convex. 
\end{thm}
\textit{Proof.} Let $v$ be as described above. Let $\epsilon>0$. Suppose there is a body $T$ whose log projection vector is very close to $v$. We may assume that it satisfies \begin{align}|T_{123}|>(1-\epsilon)\sqrt{|T_{12}||T_{23}||T_{13}|}\ \ \ \ |T_{234}|>(1-\epsilon)\sqrt{|T_{42}||T_{23}||T_{43}|}
\end{align}
Let us first show that (3.2) can not hold when $T_{123},T_{234}$ are both bounded open subsets within the respective subspaces. We now apply Lemma 10 on $T_{123}$ and $T_{234}$. We need that for every $i,j$ there is a set in the cover containing $i$ but not $j$. This is true for the above covers of $\{1,2,3\}$ and $\{2,3,4\}$ so we deduce that there are open product sets $B,B'$ in Span($\{1,2,3\})$ and Span($\{2,3,4\})$ respectively for which $$|T_{123}\triangle B|,|T_{234}\triangle B'|\le C\epsilon$$ Let $b_1,b_2,b_3,b_2',b_3',b_4'$ be the measures of the projections of $B,B'$ on the axes. Since $|B|^2=|B_{12}||B_{23}||B_{13}|$ and the $T$ we consider have $\log |T_{123}|,\log |T_{234}|$ close to the fixed numbers $v_{123},v_{234}$ we have that for some $\epsilon_1$ (which is dependent only on $\epsilon$ and as $\epsilon \rightarrow 0$ also $\epsilon_1\rightarrow 0$) $|B_{ij}\triangle T_{ij}|\le \epsilon_1$ for each $i,j\in \{1,2,3\}$. This comes from the fact that $L=T_{123}\cap B_{123}$ has volume very close to $T_{123}$ and $B_{123}$ and $L_{ij}\subset T_{ij},B_{ij}$ so we can combine that with $(6)$, the equality for $B$ and $|L|^2\le |L_{12}||L_{23}|L_{13}|$ to get $\epsilon_1$.  We may assume a similar fact holds for $B'$ by increasing $\epsilon_1$. But now we have 
\begin{align*}
\log b_2\sim v_{123}-v_{13}=-1,\log b_3\sim v_{123}-v_{12}=1\\ \log b_2'\sim v_{234}-v_{34}=1,\log b_3'\sim v_{234}-v_{24}=-1
\end{align*}
but this contradicts $|B_{23}\triangle B_{23}'|\le 2\epsilon_1$.\\
\\ 
Since $T_{123}$ is not an open subset we may approximate it by a bounded open subset above by regularity of measure. So we can have some $X$ such that $X$ is an open subset of Span ($\{1,2,3\})$ and $|X\triangle T_{123}|\le \epsilon$ Now we can approximate $X$ by Lemma 10 which then approximates $T_{123}$. Similarly we can approximate $T_{234}$ and repeat the above proof to get a contradiction. \\
\\
 The above gives us a neighbourhood of $v$ which does not contain points in $\psi_n$ but since $v\in BT_n$ and $\overline{\conv}(\psi_n)=BT_n$ we know $\psi_n\neq \conv (\psi_n)$ and hence $\psi_n$ is not convex. \QED
 \\
 \\
 We now move on to our last result. This is a lemma that says a little about the structure of $\psi_n$ itself. If we have $v\in BT_n$ and take a large $\lambda$, then the non-linearized uniform cover inequalities for $\lambda v$ still hold but the difference in magnitude between the right-hand side and left-hand side is increasing (e.g if $|T_1||T_2|>|T_{12}|$ then for large $R>1$ $|T_1|^R|T_2|^R>>|T_{12}|^R$). This can suggest that it should be easier to construct a body with log-projection vector $\lambda v$ for large $\lambda$ then one with $v$. Of course this only applies to uniform cover inequalities which are strict inequalities. This leads us to the following lemma.
 
\begin{lem}
Let $v\in BT_n$ with all of the irreducible uniform cover inequalities on $v$ holding with strict inequality. Then for all large enough $\lambda>0$ we have $\lambda v\in \psi_n$.
\end{lem}
\textit{Proof.} Let $x_A=e^{\lambda v_A}$. We want a body with projection volumes $x_A$. We know that all of the finitely many irreducible cover inequalities hold with strict inequality hence for any $C>0$ there is some $\lambda$ such that for any $k$ and any irreducible $k$-uniform cover $A_1,..,A_r$ of any $A\subset [n]$ \begin{align}\prod\limits_i x_{A_i}\ge Cx_A^k
\end{align}
Now we prove the following statement:\\
\\
\textbf{Claim.} Suppose that $S_1,S_2,...,S_r$ was a sequence of subsets of $[n]$ such that if $i<j$ then $S_j\not \subset S_i$. Then there is some $C>0$ such that if $(7)$ holds for any irreducible cover of any subset of any of the $S_i$ then there is a compact $T= \cup_{i=1}^r X_i$ (possibly with zero volume) where $X_i$ are boxes in Span$(S_i)$ all of whose projections on any subspace are disjoint, such that $|T_{S_i}|=x_{S_i}$ for all $i$.\\
\\
\textit{Proof of Claim.} We use induction on $r$. For $r=1$ the result is trivial with $C=C_1=1$. Suppose now the statement is true for the sequence $S_1,..,S_{r-1}$ with constant $C=C_{r-1}$. We now do the induction step. Since there are finitely many irreducible uniform cover inequalities. There is a constant $C=C_r$ such that if $x_{S_i}$ satisfies $(7)$ for any irreducible cover of any subset of any $S_i$ then the vector $y$ defined by \begin{align*}
    y_A=x_A/2, A\neq S_r\\ 
    y_{S_r}=x_{S_r}\ \ \ \ \ \ \ \ \ \ \ \ 
\end{align*}
satisfies all the irreducible uniform cover inequalities of any subset of $S_r$, while for any vector $z\le y$ (meaning that $z_A\le y_A$ for all $A$) $x-z$ satisfies all the inequalities in $(7)$ when $A$ is any subset of any of the sets $S_1,..,S_{r-1}$ with $C=C_{r-1}$. Now we proceed as in the proof of the Bollobás-Thomason Box Theorem [2, Theorem 1]. Let $z$ be a variable vector and consider the set of inequalities 
\begin{flalign*}
&(\textrm{i}) \ 0\le z_Y\le y_Y \textrm{ for all }Y\subset S_r &&\\
&(\textrm{ii}) \ z_Y\le \prod_{i\in Y} z_i \textrm{ for all }Y\subset S_r&&\\
&(\textrm{iii})\ y_{S_r}^k\le \prod\limits_{i=1}^s z_{Y_i} \textrm{ for each }k\textrm{-uniform irreducible cover }Y_1,..,Y_s\textrm{ of }S_r
\end{flalign*}
We see that $z=y$ is a solution and since this is a finite set of inequalities by compactness we may take a minimal solution $z$. Define $z_Y=0$ for all $Y\not \subset S_r$. Just like in the proof of Theorem 1 in [2] we have that $x_{S_r}=z_{S_r}=\prod_{i\in S_r} z_i$ and $z_Y=\prod_{i\in Y}z_i$ for all $Y\subset S_r$. Now we let $X_r$ be a box in Span($S_r$) with sidelengths $z_i$. By induction there is a body $T'=\cup_{i=1}^{r-1}X_i$ which is a union of boxes $X_i$ in Span($S_i$) all of whose projections are disjoint and such that for $1\le i \le r-1$ $|T'_{S_i}|=(x-z)_{S_i}$. Thus we can construct $T=T'\cup X_r$ such that we put $X_r$  far away from $T'$ in all directions so that no projections of $X_r$ intersect any projections of $T'$. This will insure that $x_{S_i}=T_{S_i}$ for all $i$. This proves our claim.\\
\\
Now we consider the non-empty subsets $S_1,..,S_{2^n-1}$ of $[n]$ ordered by size (arbitrary for sets of the same size). This ordering satisfies the condition of the claim and hence we can pick a $\lambda$ for which $x$ satisfies all the inequalities of type $(7)$ for $C=C_{2^n-1}$ given by the claim. Then by the claim there is a compact $T$ such that $T$ is a body with log projection vector $\lambda v$.\QED\\
\\
Surprisingly, this lemma gives another proof of Theorem 4.\\
\\
\textit{Second proof of Theorem 4.} Let $v\in BT_n$ and let $\epsilon>0$. The vector $w$ defined by $w_A=v_A+\epsilon$ for all $A$ has the property that all irreducible uniform cover inequalities hold with a strict inequality because $v\in BT_n$ and any such inequality has more terms on the larger side of the inequality. Thus by Lemma 12 there is a $\lambda$ such that $\lambda w\in \psi_n$, but $\bf{0}\in \psi_n$ implies $w\in \conv (\psi_n)$. Since $\epsilon$ was arbitrary, $v\in \overline{\conv}(\psi_n)$.\QED\\
\\
From the proof of Theorem $11$ we have a vector $v\in BT_n$ with a neighbourhood around it disjoint from $\psi_n$. By the above there is some $w$ in that neighbourhood which has all irreducible uniform cover inequalities strict, meaning some positive multiple of $w$ is in $\psi_n$. This tells us that $\psi_n$ is not a cone. What else can we say about $\psi_n$? We shall end by stating the following $2$ conjectures we have not resolved, the second one of which naturally arises from Lemma $12$.\\
\\
We close with two
conjectures. The first is that the set $\psi_n$ is closed. It is very frustrating that we cannot answer this question. One would hope for some kind of
compactness argument. Now, if there were some kind of `canonical form' for a
body that realises a given point in $\psi_n$ then this might be achievable, but
unfortunately the work of Tan and Zeng [6] on `rectangular flowers' makes it unlikely that there is such a canonical form, at least for $n \geq 4$.
\begin{conj}
The set $\psi_n$ is closed.
\end{conj}
The second conjecture is a natural strengthening  of Lemma 12. It seems that scaling up a vector should only make it `more constructible'.
\begin{conj}
If $v\in \psi_n$ and $\lambda > 1$ then $\lambda v\in \psi_n$.
\end{conj}

\section{References}
[1] Bollobás, B. \textit{Linear Analysis: An Introductory Course,} Cambridge University Press, 1990.
\newline
\newline
 [2] B. Bollobás and A. Thomason, Projections of bodies and hereditary properties of hypergraphs, Bull. London Math. Soc. \textbf{27}(1995), 417–424.\newline
 \newline
[3] F. R. K. Chung, R. L. Graham, P. Frankl and J. B. Shearer. Some intersection theorems for ordered
sets and graphs. J. Combinatorial Theory (A) \textbf{43}(1986), 23–37.
\newline
\newline 
[4] D. Ellis, E. Friedgut, G. Kindler and A. Yehudayoff, Geometric stability via information theory, arXiv:1510.00258.\newline
\newline
[5] D. Luenberger, \textit{Introduction to Linear and Non-Linear Programming},
Addison-Wesley 1984.
\newline
\newline
[6] Z. Tan and L. Zeng, On the inequalities of projected volumes and the constructible
region,  	arXiv:1410.8663.
\end{document}